\documentclass{amsart}
\usepackage{epsfig,amsmath,amsthm,latexsym,graphicx,amsfonts,amssymb}
\usepackage[OT2,OT1]{fontenc}
\newcommand\cyr{%
\renewcommand\rmdefault{wncyr}%
\renewcommand\sfdefault{wncyss}%
\renewcommand\encodingdefault{OT2}%
\normalfont
\selectfont}
\DeclareTextFontCommand{\textcyr}{\cyr}
\newtheorem{theorem}{Theorem}[section]
\newtheorem{lemma}[theorem]{Lemma}
\newtheorem{corollary}[theorem]{Corollary}
\newtheorem{proposition}[theorem]{Proposition}

\newcommand{\Vol}{\operatorname{Vol}}
\newcommand{\BH}{\mathbb{H}}
\newcommand{\BR}{\mathbb{R}}
\newcommand{\BC}{\mathbb{C}}

\newcommand{\mr}[1]{\mathrm{#1}}

\begin{document}
\title{Minimum-volume hyperbolic 3-manifolds}
\author{Peter Milley}\thanks{Partially supported by NSF grant
DMS-0554624 and by ARC Discovery grant DP0663399.}\address{Department
of Mathematics and Statistics\\University of Melbourne\\Melbourne,
Australia}\email{P.Milley@ms.unimelb.edu.au}

\maketitle

\section{Introduction}
The classification of small-volume hyperbolic 3-manifolds has been an
active problem for many years, ever since Thurston suggested that
volume was a measure of the complexity of a hyperbolic
3-manifold. Quite recently in \cite{gmm3} the author along with David
Gabai and Robert Meyerhoff used a geometrical construction called a
\emph{Mom-$n$ structure} to tackle the classification problem, and
succeeded in showing the following:

\begin{theorem}(Gabai, Meyerhoff, and Milley)
\begin{itemize}
\item If $M$ is an orientable one-cusped hyperbolic 3-manifold with
volume less than $2.848$, then $M$ can be obtained by a Dehn filling
on one of the following manifolds in the SnapPea census: m125, m129,
m202, m203, m292, m295, m328, m329, m359, m366, m367, m391, m412,
s596, s647, s774, s776, s780, s785, s898, or s959.
\item If $M$ is a closed orientable hyperbolic 3-manifold with volume
less than that of the Weeks manifold, then $M$ can be obtained by a
Dehn filling on a one-cusped manifold with volume less than $2.848$,
and hence can be obtained by a Dehn filling on one of the 21 manifolds
listed above.
\end{itemize}\qed
\label{thrm:gmm3}
\end{theorem}
Note that the Weeks manifold is the manifold obtained by $(2,1)$
Dehn filling on the manifold m003 in the SnapPea census. We should
also point out that strictly speaking the list of manifolds in the
first part of Theorem \ref{thrm:gmm3} could be shortened considerably,
as the manifolds on the list from m125 through m391 are all Dehn
fillings on the manifold s776. However, we present the list in full to
emphasize that the 21 manifolds in the theorem are exactly the
Mom-$2$'s and Mom-$3$'s defined in \cite{gmm2}, and that another way
of stating Theorem \ref{thrm:gmm3} is to say that a suitably small
one-cusped or closed hyperbolic 3-manifold must possess a Mom-$2$ or
Mom-$3$ structure.

This theorem raises an obvious question: exactly which small-volume
manifolds can be obtained by a Dehn filling on one of these 21
manifolds? The purpose of this paper is to answer that question, and
in so doing to complete the proof of the following long-standing
conjectures:

\begin{theorem}
If $N$ is a complete orientable one-cusped hyperbolic 3-manifold with
$\Vol(N)\le 2.848$ then $N$ is homeomorphic to one of the following
manifolds in the SnapPea census: m003, m004, m006, m007, m009, m010,
m011, m015, m016, or m017.
\label{thrm:one_cusped}
\end{theorem}

\begin{theorem}
If $N$ is a closed orientable hyperbolic 3-manifold with $\Vol(N)\le
0.943$ then $N$ is homeomorphic to the Weeks manifold.
\label{thrm:compact}
\end{theorem}

The proof of the above two theorems breaks into several steps. In
section 2 we show how we can restrict our attention to a finite
number of Dehn fillings on any given manifold, and use the computer
program Snap (\cite{g}) to make a preliminary estimate of how many of
those fillings result in hyperbolic 3-manifolds. In section 3 we use
the work of Harriet Moser to rigorously confirm hyperbolicity for
those manifolds, and also to confirm which of those hyperbolic
manifolds have suitably small volumes. Finally in section 4 we confirm
that there are no additional hyperbolic manifolds to consider other
than the ones found by Snap.

A word about availability: the results in this paper were proved with
extensive computer assistance, utilizing several different
programs, both pre-existing and written by the author. These include
PARI/GP (\cite{coh}), SnapPea (\cite{w}), Snap (\cite{g}), a custom
PARI/GP routine written by Harriet
Moser and modified by the author, programs in C++ and Perl written
by author, and one or two Unix shell scripts. Rather than
include all of the intermediate results from all of these varied
programs, instead we only include detailed results from the initial
step of the analysis regarding the manifold s776 in the next section;
in all other cases the results of various computations are described
rather than given in detail. Full details of the computations along
with the various original programs used to produce them are available
from the author upon request (\cite{milley}).

The author wishes to express his thanks to Robert Meyerhoff for his
invaluable assistance with these results, and to Craig Hodgson for his
valuable comments.

\section{Dehn filling bounds}

The following result was recently proved in \cite{fkp}:
\begin{theorem}(Futer, Kalfagianni, and Purcell) Let $M$ be a
complete, finite-volume hyperbolic manifold with cusps. Suppose $C_1$,
\ldots, $C_k$ are disjoint horoball neighborhoods of some subset of
the cusps. Let $s_1$, \ldots, $s_k$ be slopes on $\partial C_1$,
\ldots, $\partial C_k$, each with length greater than $2 \pi$. Denote
the minimal slope length by $l_{\mr{min}}$. If $M(s_1,\ldots,s_k)$
satisfies the geometrization conjecture, then it is a hyperbolic
manifold, and
\[
\Vol(M(s_1,\ldots,s_k)) \ge \left(1-\left(\frac{2\pi}
{l_{\mr{min}}}\right)^2\right)^{3/2} \Vol(M).
\]
\label{thrm:fkp}
\end{theorem}

We wish to use this theorem to find a bound on $l_{\mr{min}}$ in the
case where the filled manifold has small volume, hence we re-arrange
the above result as follows:
\begin{corollary}Suppose $M$, $s_1$, \ldots, $s_k$, and $l_{\mr{min}}$
are defined as above, and let $N=M(s_1,\ldots,s_k)$. Assuming the
geometrization conjecture, if $N$ is hyperbolic we have
\[
l_{\mr{min}} \le
2\pi \left( \sqrt{1- \left(\frac{\Vol(N)}{\Vol(M)}\right)^{2/3}}
\right) ^{-1}\,.
\]
\label{cor:fkp}
\end{corollary}

We begin our analysis with the ``magic manifold'' s776, which is the
complement in $S^3$ of the alternating three-element chain link. This
manifold has already been extensively analyzed in \cite{mp}, and we
will refer to some of the results of that paper in what follows. This
is also the only three-cusped manifold of the ones listed in Theorem
\ref{thrm:gmm3}; the remainder of those manifolds have two cusps. If
we fill in two of the cusps of s776 we may obtain a one-cusped $N$
with $\Vol(N)\le 2.848$; our first goal is to determine what one of
the two Dehn filling coefficients must be in that case.

Suppose that $N$ is a one-cusped
hyperbolic 3-manifold obtained by Dehn filling on two of the three
cusps of s776, and suppose that $\Vol(N)\le 2.848$. Since
$\Vol(\mr{s}776)=5.333\ldots$, Corollary \ref{cor:fkp} implies that
the minimum slope length $l_{\mr{min}}$ must be less than or equal to
$10.746\ldots\ $. Let $C_1$ and $C_2$ be horoball neighbourhoods of
the two filled cusps, with equal volume but chosen to be large as
possible while still having disjoint interiors. Then coordinates can
be chosen on the boundary of each $C_i$ such that the slope of the meridian
corresponds to the complex number $(1+i\sqrt{7})/2$ while the
slope of the longitude corresponds to $2$. This implies that, up to
symmetry around the origin, the Dehn
filling coefficients on one of the two filled cusps must one of the
pairs listed in the table in
\begin{figure}
\begin{center}
\begin{tabular}{|c|c|c|c|c|c|}
\hline
$(-8,1)$ & $(-8,3)$ & $(-7,1)$ & $(-7,2)$ & $(-7,3)$ & $(-7,4)$ \\
$(-6,1)$ & $(-6,5)$ & $(-5,1)$ & $(-5,2)$ & $(-5,3)$ & $(-5,4)$ \\
$(-4,1)$ & $(-4,3)$ & $(-4,5)$ & $(-3,1)$ & $(-3,2)$ & $(-3,4)$ \\
$(-3,5)$ & $(-2,1)$ & $(-2,3)$ & $(-2,5)$ & $(-1,1)$ & $(-1,2)$ \\
$(-1,3)$ & $(-1,4)$ & $(-1,5)$ & $(0,1)$ & $(1,0)$ & $(1,1)$ \\
$(1,2)$ & $(1,3)$ & $(1,4)$ & $(1,5)$ & $(2,1)$ & $(2,3)$ \\
$(3,1)$ & $(3,2)$ & $(3,4)$ & $(4,1)$ & $(4,3)$ & $(5,1)$ \\
$(5,2)$ & $(5,3)$ & $(6,1)$ & $(7,1)$ & & \\
\hline
\end{tabular}
\end{center}
\caption{The possible Dehn filling coefficients on one of the two
filled cusps of s776.}
\end{figure}
figure 1.

Consulting \cite{mp}, we see that the fillings whose coefficients in
these coordinates are $(-2,1)$, $(-1,1)$ $(0,1)$, $(1,0)$, and $(1,1)$
will result in non-hyperbolic manifolds, and further that there is no
way to fill the remaining cusps of these five non-hyperbolic manifolds
to get a hyperbolic result. The remaining 41 sets of Dehn
filling coefficients, if used to fill in only one cusp of s776, will result in
hyperbolic two-cusped manifolds.

(Note that coordinate system used to describe the filling coefficients
here is SnapPea's, and is based on a choice of meridian and longitude
on each cusp where the meridian is a shortest homotopically
non-trivial curve on the cusp torus and the longitude is a shortest
homotopically non-trivial curve which is linearly independent from the
meridian. Martelli and Petronio use a different coordinate system to
describe Dehn fillings in \cite{mp}, corresponding to a presentation
of s776 as the complement of a link in $S^3$, and it is necessary to
translate between the two systems to get the above result.)

Using Snap, it is possible to identify all but 4 of the hyperbolic
manifolds described above as coming from the SnapPea census (although
in a few cases it is necessary to retriangulate the manifold
first). The four exceptions correspond to the coefficients $(-8,1)$,
$(-7,1)$, $(6,1)$, and $(7,1)$. These four filling coefficients
produce four non-isometric two-cusped hyperbolic manifolds which are
not in the census, presumably due to high Matveev complexity. Of the
manifolds that do appear on the census, 13 already appear on the list
of 2-cusped Mom-2's and Mom-3's. The remaining 24 manifolds are listed
in
\begin{figure}
\begin{center}
\begin{tabular}{|c|c|c|c|c|c|c|c|}
\hline
m357 & m388 & s441 & s443 & s503 & s506 & s548 & s549 \\
s568 & s569 & s576 & s577 & s578 & s579 & s601 & s602 \\
s621 & s622 & v1060 & v1061 & v1178 & v1180 & v1203 & v1204 \\
\hline
\end{tabular}
\end{center}
\caption{The 24 2-cusped manifolds resulting from relevant surgeries
on s776 which are not Mom-$2$'s or Mom-$3$'s.}
\end{figure}
figure 2.

Next we wish to analyze the Dehn surgery spaces of the 20
2-cusped manifolds in the list of Mom-$2$'s and Mom-$3$'s, the 24
manifolds listed in figure 2 which were obtained by filling in one
cusp of the manifold s776, and the four unnamed manifolds obtained by
filling in one cusp of s776 with Dehn filling coefficients $(-8,1)$,
$(-7,1)$, $(6,1)$, or $(7,1)$. For each manifold our goal is to
rigorously enumerate all possible fillings on one cusp which
result in a 1-cusped hyperbolic 3-manifold $N$ with $\Vol(N)\le
2.848$. With one exception, each of these 48 manifolds admits a
symmetry which exchanges the two cusps, and hence it does not matter
which cusp we choose to fill; the one exception is the manifold
s785, for which the analysis that follows has to be performed twice,
once for each cusp.

For each 2-cusped manifold $M$ the procedure for performing the
analysis is the same. First we use Corollary \ref{cor:fkp} to obtain
an upper bound on the length of the slope that needs to be considered,
then use this upper bound to obtain a finite list of possible Dehn
filling coefficients, just as in the previous section. This requires
three pieces of data which we obtain from SnapPea: the volume of $M$
and the slopes corresponding to a meridian and longitude on a maximal
horoball neighbourhood of the cusp being filled. While this data
cannot be considered exact, at this stage of the process we are only
concerned with ensuring that no relevant Dehn filling coefficients are
omitted; hence ``fudging'' upward the value of $l_{\mr{min}}$ obtained
from Corollary \ref{cor:fkp} should compensate for any floating-point
error in this step. This process is readily automated and results in
1278 candidate Dehn fillings to consider.

Having established a finite list of one-cusped manifolds to consider,
we now move on to the closed case. We will anticipate the solution to
Theorem \ref{thrm:one_cusped} and only examine those manifolds that
can be obtained by Dehn filling on one of the first 10 orientable
one-cusped hyperbolic 3-manifolds in the SnapPea census. The procedure
is identical to that used above: for each of the 10 one-cusped
manifolds, use Corollary \ref{cor:fkp} to obtain an upper bound on the
length of the surgery slopes and consequently on the number of Dehn
surgery coefficients to be considered. This results in 224 total Dehn
fillings that need to be examined.

In both cases, the Dehn fillings need to be classified into fillings
which result in hyperbolic manifolds and fillings which result in
non-hyperbolic manifolds; the hyperbolic manifolds then have to be
further classified by volume to determine which of the one-cusped
(resp. closed) filled manifolds have volume less than or equal to
$2.848$ (resp. $0.943$). The computer program Snap can perform this
classification for us, but using floating-point arithmetic; for the
sake of rigour, Snap's results need to be confirmed. This is the goal
of the next two sections of this paper. For now we simply state Snap's
results.

Of the $1278$ Dehn fillings on one-cusped manifolds, Snap claims to
find positively oriented ideal triangulations on $989$ of them and
negatively oriented or partially flat ideal triangulations on another
$53$. (The distinction between these types of triangulations will
become relevant when we confirm hyperbolicity in the next section.)
Snap fails to find a hyperbolic ideal triangulation on the remaining
$236$ manifolds. Of the $224$ Dehn fillings on one-cusped manifolds,
Snap claims to find positively oriented ideal triangulations on $137$
of them, negatively oriented or partially flat ideal triangulations on
$18$, and fails to find a hyperbolic ideal triangulation on the
remaining $69$.

\section{Confirming hyperbolicity and volume}

Snap attemps to find a hyperbolic ideal triangulation on a filled
manifold by using Newton's Method to solve a system of gluing
equations that arise from studying the holonomy around the edges of an
ideal triangulation (for details of this process see for example
\cite{mos} or \cite{thu}). Consequently for each of the filled
manifolds for which Snap claims to have found a solution there is a
legitimate question as to whether the solution truly exists, or
instead is a result of floating-point error.  While Snap can attempt
to find an exact algebraic description of a hyperbolic manifold, there
is no guarantee that it will be able to so for every manifold that we
are interested in. Fortunately Harriet Moser has created an algorithm
(\cite{mos}) which can confirm whether or not a genuine solution to
the gluing equations exists in a small neighbourhood of the solution
calculated by Snap using Newton's Method. This algorithm takes the
gluing equations and Snap's computed solution and uses this data to
evaluate a pair of inequalities; either inequality, if true, confirms
the existence of a true solution.

For our purposes the only technical shortcoming of Moser's algorithm
is that it requires that the hyperbolic structure computed by Snap
uses only positively oriented tetrahedra. Several of the filled
manifolds that we wish to prove hyperbolic have Snap-computed
hyperbolic ideal triangulations which include negatively oriented or
flat tetrahedra. Fortunately in all but five cases, these manifolds
were shown by Snap to be homeomorphic to manifolds with positively
oriented triangulations. (The algorithms Snap uses to do this operate
on the combinatorial data incorporated in ideal triangulations and
therefore are not subject to floating-point inaccuracy.)  The five
exceptions are all closed filled manifolds; all of the $53$ one-cusped
manifolds with negatively oriented or partially flat triangulations
were successfully re-triangulated, resulting in $39$ new positively
oriented triangulations. Note that some duplicates were eliminated at
this stage. Of the $18$ closed filled manifolds with negatively
oriented or partially flat triangulations, $13$ were successfully
re-triangulated. Of the five exceptions, three of them in fact turned
out not to be hyperbolic. These three manifolds are m004(1,0),
m009(1,0) and m015(1,0). Snap classifies these manifolds as
``nongeometric'', and finds shape parameters for their ideal
triangulations which seem to have been obtained by analytic
continuation along some path which is not contained in the space of
hyperbolic Dehn fillings. A cursory examination of m004, m009, and
m015 shows that (1,0) surgery on these manifolds results in a
nonhyperbolic closed manifold, so these three filled manifolds were
added to the list of non-hyperbolic manifolds which are dealt with in
the next section.

The remaining two exceptions were both homeomorphic to the manifold
known as Vol3, or m007(3,1) in the SnapPea census. All known ideal
triangulations of this manifold have negatively oriented tetrahedra,
and it is an open problem whether or not this is true for all ideal
triangulations of this manifold. Since our goal is only to prove that
there are no closed hyperbolic manifolds smaller than the Weeks
manifold, we will satisfy ourselves with just proving that Vol3 is not
smaller than the Weeks manifold when we address the question of volume
shortly. As for those manifolds for which positively oriented ideal
triangulations can be found, Moser's algorithm was applied to all of
them and in all cases the algorithm confirmed the hyperbolicity of the
manifold.

We now turn our attention to the problem of rigorously computing the
volume of a hyperbolic 3-manifold. For a manifold with an ideal
triangulation with positively oriented tetrahedra, the volume is given
by the well-known formula (see, for example, \cite{mil}):
\[
\Vol(M) = \sum_\tau
\textcyr{L}(\theta_{\tau,1}) + \textcyr{L}(\theta_{\tau,2}) +
\textcyr{L}(\theta_{\tau,3})
\]
where $\tau$ varies over all simplices in the triangulation, and
$\theta_{\tau,i}$ are the dihedral angles associated with $\tau$ for
$i=1$ to $3$. If $z$ is the complex parameter corresponding to $\tau$,
then the three dihedral angles associated to $\tau$ are the arguments
of the complex numbers $z$, $1/(1-z)$, and $(z-1)/z$. The function
$\textcyr{L}(\theta)$ is the \emph{Lobachevsky function}:
\[
\textcyr{L}(\theta) = \int_0^\theta -\log |2 \sin t|\, dt\, .
\]
While Snap will compute this function for us, this computation uses
floating-point arithmetic and hence cannot be considered completely
accurate. Moreover, Snap's computations are based on its own solution
to the gluing equation for the manifold, which may not be accurate to
begin with. What we require is an upper bound on the difference
between Snap's computed volume and the actual volume of the manifold,
so that we can determine a rigorous lower bound on the volume of the
manifold.

\cite{mos} provides the first half of the solution to this
problem. One of the intermediate steps in Moser's algorithm produces
an upper bound $\delta$ on the distance in $\BC^n$ between Snap's
computed solution to the filling equations and the actual
solution. (While this number varies from manifold to manifold it is
usually on the order of $10^{-3}$ or smaller.) This number represents
the error associated to the input to the volume function, and it is
trivial to modify Moser's algorithm to produce this number for a given
manifold.

It then remains to compute a range of possible values of $\Vol(M)$
given the quantity $\delta$ and the tetrahedron shapes computed by
Snap. For this we turn to \emph{affine 1-jets}, as described in
\cite{gmm3} and \cite{gmt}.

We refer the reader to those papers for details, but a brief
description of the concept is as follows. An affine 1-jet is a linear
function $j:[-1,1]^n\rightarrow\BR$ together with an error term,
usually described as a tuple $(j_0;j_1,\ldots,j_n;j_\epsilon)$ where
$j_0$ is the constant coefficient of the function, $j_1$, \ldots,
$j_n$ are the linear coefficients, and $j_\epsilon$ is the error
term. (In \cite{gmt} complex affine 1-jets were used, where the domain
of the function described is $U^n$ where $U$ is the unit disk
$\{|z|\le 1\}$ in $\BC$, but the principle is the same.)
Mathematically, an affine 1-jet represents a neighbourhood of the
linear function $(x_1,\ldots,x_n) \mapsto j_0+\sum j_i x_i$ in the
space of functions equipped with the sup norm. From a practical
perspective affine 1-jets are well suited to use as approximations to
non-linear functions since it is possible to construct definitions for
the basic arithmetic operations such that, for example, the affine
1-jet corresponding to $j+k$ contains all sums $f+g$ where $f$ and $g$
are any functions in the neighbourhoods corresponding to $j$ and $k$
respectively. Using facts about the IEEE standard for floating-point
arithmetic, it is even possible to define these operations in such a
way that they take into account the floating-point error that may
occur during the computation. Hence using affine 1-jets it is possible
to program a computer to compute the value of any rational polynomial
and keep track of its own error, allowing for rigorous computer-aided
proofs of inequalities involving rational polynomials. In addition,
using Taylor approximations and Taylor's theorem allows this technique
to be extended to functions whose derivatives are rational
polynomials, such as was done with the logarithm function in
\cite{gmm3}.

Clearly the formula for $\Vol(M)$ given above is not a rational
polynomial, nevertheless we can adapt it to be suitable for
computation with affine 1-jets. The first difficulty is computing the
dihedral angles $\theta_{z,i}$ from $z$ for $i=1$, $2$, $3$, which
requires the arctangent function. Fortunately the derivatives of the
arctangent function are all rational polynomials, so a Taylor
approximation can be applied in the same manor as for the logarithm
function. The second difficulty is computing the Lobachevsky function
itself. For that, we use the following well-known series expansion
(see, for example, \cite{mil}):
\begin{lemma}
\[
\textcyr{L}(\theta) = \theta\left(
1 - \log|2\theta| + \sum_{n=1}^\infty
\frac{|B_{2n}| (2\theta)^{2n}}{2n(2n+1)!}
\right)
\]
where $\{B_{2n}\}$ are the Bernoulli numbers.
\end{lemma}
(Note that Milnor uses a slightly older definition of the Bernoulli
numbers and consequently Milnor's formula uses the notation $B_n$
instead of $|B_{2n}|$; nevertheless the meaning is the same.)

Since the logarithm function for affine 1-jets was implemented in
\cite{gmm3} we can use a truncated version of this series expansion to
calculate $\textcyr{L}(\theta)$, but since the derivatives of
$\textcyr{L}(\theta)$ are not rational polynomials we can't use
Taylor's theorem to get a bound on the error. Instead, we make use of
the following lemma:

\begin{lemma}
Let
\[
l_n = \frac{|B_{2n}| 2^{2n}}{2n(2n+1)!}
\]
for $n=1$, $2$, \ldots\ . Then $l_n/l_{n+1}>\pi^2$ for all $n\ge 1$.
\end{lemma}
\begin{corollary}
If $\theta^2<\pi^2/2$ then
\[
\sum_{n=k}^\infty l_n \theta^{2n} < 2 l_k \theta^{2k}
\]
for all $k\ge 1$.
\end{corollary}

\noindent\emph{Proof:} To prove the lemma, note the following fact
about the Bernoulli numbers (see e.g. \cite{woo}):
\[
B_{2n} = \frac{(-1)^{n-1} 2 (2n)!}{(2\pi)^{2n}} \zeta(2n)
\]
where $\zeta(x)$ is the Riemann zeta function. Substituting this into
the definition of $l_n$ we get
\[
l_n = \frac{\zeta(2n)}{\pi^{2n}n(2n+1)}
\]
and hence
\[
\frac{l_n}{l_{n+1}} = \pi^2 \frac{(n+1)(2n+3)\zeta(2n)}
{n(2n+1)\zeta(2n+2)}\, .
\]
Since $\zeta(x)$ is a decreasing function for $x>1$, the right-hand
side of the above equation is clearly greater than $\pi^2$, proving
the lemma. The lemma further implies that $l_n\le l_k/\pi^{2(n-k)}$
for all $n\ge k$; the corollary follows.\qed

\smallskip
The above corollary allows us to compute a bound on the error term
when using affine 1-jets to evaluate the series expansion for
$\textcyr{L}(\theta)$. Specifically, if $|\theta|$ is
provably smaller than $\pi/\sqrt{2}$, and if any term of the series
expansion can be proved to be less than $\epsilon$ in magnitude, where
$\epsilon$ is any arbitrarily small positive constant, then we know
the that sum of all remaining terms of the series must be at most
$2\epsilon$.

Since $\textcyr{L}(\theta)$ is both odd and $\pi$-periodic, in theory
for any real number $\theta$ we can always find a new argument
$\theta_0$ such that $0\le\theta_0\le\pi/2 < \pi/\sqrt{2}$ and
$\textcyr{L}(\theta)=\textcyr{L}(\theta_0)$. For affine 1-jets,
however, this transformation may not be possible if the error
associated to the original 1-jet is large. Nevertheless our program
used such a transformation to evaluate $\textcyr{L}(\theta)$ for
1-jets and was programmed to report an error if the condition
$|\theta_0|<\pi/\sqrt{2}$ could not be guaranteed. Fortunately, this
error condition never arose.

Using these techniques, a total of $1028$ one-cusped hyperbolic
3-manifolds with positively oriented ideal triangulations were
analyzed; of these, all but $48$ manifolds were proved to have volume
greater that $2.848$. Of the $48$ exceptions, $47$ can be proved by
Snap to be isometric to one of the 10 cusped census manifolds with
volume less than $2.848$, i.e. one of the 10 manifolds listed in
Theorem \ref{thrm:one_cusped}. The remaining manifold turns out to be
isomorphic to m019, which according to Snap should have volume greater
than $2.944$. This was the one cusped case where the error term
calculated by our program was large enough to prevent it from making a
correct determination. However, applying our program to a different
triangulation of m019 proved that this manifold did indeed have volume
greater than $2.848$, completing the results of this section for the
one-cusped filled manifolds.

For the closed case, a total of 150 closed filled manifolds were
analyzed: 137 manifolds for which Snap found a positively oriented
triangulation originally, and another 13 for which a positively
oriented re-triangulation could be found. All but eight of the
triangulations were proved to have volume greater than $0.943$. Three
of the exceptions were isometric to the Weeks manifold, as
expected. The remaining five exceptions were isometric to Vol2, the
second smallest known closed orientable hyperbolic 3-manifold, which
according to Snap has volume greater than $0.981$. As in the cusped
case the error term calculated by our program was large enough to
prevent it from making a correct determination, but re-applying our
program to a different triangulation of Vol2 proved that this manifold
has volume greater than $0.943$.

Finally, we need to prove that the manifold Vol3 has volume greater
than $0.943$ as discussed previously. (For a detailed proof that Vol3
is hyperbolic, see for example \cite{JR}.) Although this manifold does
not have any known positively oriented ideal triangulation, with some
work it is possible to persuade Snap to find a positively oriented
ideal triangulation of its unique double cover, shown by Snap to be
isometric to m036(-3,2). Applying Moser's algorithm to confirm the
hyperbolicity of this cover and then applying our rigorous program to
compute volume shows that the volume of the double cover of Vol3 is
demostrably greater than $1.886$, completing our results for this
section.

\section{Confirming non-hyperbolicity}

For the manifolds which Snap fails to find a hyperbolic structure (and
for three of the closed filled manifolds described in the previous
section) it is again necessary to confirm Snap's results. The author
knows of no automated tool comparable to Moser's algorithm to prove
the non-hyperbolicity of a manifold. However in practice the manifolds
in question can be proved to by non-hyperbolic by examining their
fundamental groups.

Specifically, in each case if we assume that $M$ is hyperbolic then we
can get a contradiction. For if $M$ is a closed or cusped hyperbolic
manifold then $\pi_1(M)$ must be a discrete finite-covolume subgroup of
$\mr{PSL}(2,\BC)$ containing no elliptic elements; that is, $\pi_1(M)$
contains no non-trivial elements which fix any points in $\BH^3$, or
equivalently no elements of $\mr{PSL}(2,\BC)$ whose trace is real and
lies in the interval $[-2,2]$. Such groups have very restrictive
properties, which we list below without proof. In the lemmas
below, $[x,y]$ denotes the commutator $xyx^{-1}y^{-1}$ of $x$
and $y$, and $1$ denotes the identity element of the group.

\begin{lemma}
Let $\Gamma$ be a discrete finite-covolume subgroup of
$\mr{PSL}(2,\BC)$ with no elliptic elements (e.g., the fundamental
group of a complete finite-volume 3-manifold), and let $a$ and $b$ be
distinct elements of $\Gamma$.
\begin{enumerate}
\item The centre of $\Gamma$ is trivial; in particular, $\Gamma$ is
not abelian. (Note we require the finite-covolume assumption here.)
\item The maximal abelian subgroups of $\Gamma$ have rank 1 or
2; if $H_1$, $H_2$ are two such maximal abelian subgroups then either
$H_1=H_2$ or $H_1\cap H_2=\{1\}$.
\item $[a,b]=1$, that is $a$ and $b$ commute, if and
only if $a$ and $b$ lie in the same maximal abelian subgroup.
\item If $a$ and $bab^{-1}$ both lie in the same maximal abelian
subgroup of $\Gamma$ then $[a,b]=1$, i.e. $b$ also lies in the same
maximal abelian subgroup.
\item $[a^n,b^m]=1$ for non-zero integers $m$ and $n$ if and only if
$[a,b]=1$.
\end{enumerate}
\label{lem:easy_groups}
\end{lemma}

In addition, we have the following simple lemma:

\begin{lemma}
If $a$, $b\in\Gamma$ and $a^nb^m=b^{-m}a^k$ for some integers $n$,
$m$, and $k$, then $[a^{n+k},b^m]=1$. If $\Gamma$ has the presentation
$\langle a,b\ |\ a^nb^ma^{-k}b^m\rangle$ then $\Gamma$ is not the
fundamental group of a complete finite-volume hyperbolic 3-manifold.
\label{lem:group_lem}
\end{lemma}

\noindent\emph{Proof:} We have
\begin{eqnarray*}
a^{n+k}b^ma^{-n-k} &=& a^k(a^nb^m)a^{-n-k} \\
                   &=& a^kb^{-m}a^{-n} \\
                   &=& a^k(a^nb^m)^{-1} \\
                   &=& a^ka^{-k}b^m \\
                   &=& b^m
\end{eqnarray*}
Suppose $\Gamma$ has the given presentation and is the fundamental
group of a complete finite-volume hyperbolic 3-manifold. Then
$[a^{n+k},b^m]=1$ and hence $[a,b]=1$ by the previous lemma. Hence
$\Gamma$ is abelian, a contradiction.\qed

Armed with these facts, many of the conjecturally non-hyperbolic
manifolds can quickly be proved to be non-hyperbolic after a cursory
examination of their fundamental group. For example, suppose the group
$\langle a,b\ |\ a^3b^2\rangle$ were the fundamental group of a
finite-volume complete hyperbolic 3-manifold. Then since $[a^3,b^2]$
is clearly trivial we must also have $[a,b]=1$ which would imply that
the whole group is abelian, a contradiction. Therefore any manifold
with fundamental group $\langle a,b\ |\ a^3b^2\rangle$ cannot be a
finite-volume complete hyperbolic 3-manifold. Similarly, the group
$\langle a,b\ |\ a^2b^2a^{-1}b^2\rangle$ cannot be the fundamental
group of a finite-volume hyperbolic 3-manifold because by the above
lemma with $n=2$, $m=2$, and $k=1$ we have $[a^3,b^2]=1$, leading to
the same contradiction.

These arguments are sufficient to confirm non-hyperbolicity in the
majority of cases. A small number of fundamental groups required
further analysis, as discussed below.

For the one-cusped filled manifolds, 236 of the Dehn fillings produced
in section 2 result in manifolds for which Snap fails to find a
hyperbolic structure. To
reduce the size of this list, a program was written using the SnapPea
kernel to re-triangulate these manifolds in an attempt to find
manifolds with identical triangulations. This effort reduced the list
to 79 conjecturally non-hyperbolic manifolds. In 68 cases the
presentation is of the form $\langle a,b|r\rangle$ where the single
relation $r$ is either of the form $a^nb^m$ for some $n$ and $m$ or of
the form $a^nb^ma^{-k}b^m$ for some $n$, $m$, and $k$ (possibly after
switching the roles of $a$ and $b$). These groups cannot be the
fundamental groups of a finite-volume hyperbolic 3-manifold: in the
first case the group would either have torsion or a non-trivial
centre, while in the second case Lemma \ref{lem:group_lem}
applies. The remaining 11 cases are analyzed in the following lemma:

\begin{proposition}
None of the following 11 groups is the fundamental group of a
complete finite-volume hyperbolic 3-manifold:
\begin{itemize}
\item $\langle a,b\ |\ ab^{-1}a^{-2}b^{-1}ab^2 \rangle$
\item $\langle a,b\ |\ ab^{-3}ab^2a^2b^2 \rangle$
\item $\langle a,b\ |\ ab^2ab^{-1}ab^2a^2b^2ab^{-1} \rangle$
\item $\langle a,b\ |\ a^3b^2a^3ba^{-2}b \rangle$
\item $\langle a,b,c\ |\ a^2cb^{-1}cb,\ b^2c^2 \rangle$
\item $\langle a,b,c\ |\ bc^{-1}b^{-1}c,\ ac^2ba^{-1}b \rangle$
\item $\langle a,b\ |\ ab^2a^3b^2a^3b^2a^3b^2 \rangle$
\item $\langle a,b\ |\ a^2b^2a^2b^2a^3b^2a^2b^2a^2b^2ab^{-1}ab^2 \rangle$
\item $\langle a,b\ |\ a^2ba^2b^{-2}a^2b \rangle$
\item $\langle a,b\ |\ ab^{-3}ab^4a^2b^4 \rangle$
\item $\langle a,b,c\ |\ a^2b^{-1}c^{-1}a^{-1}bc,\ b^2c^2 \rangle$
\end{itemize}
\label{prop:cusped_groups}
\end{proposition}

\noindent\emph{Proof:} For each of the 11 groups we assume that
the group is the fundamental group $\Gamma$ of a complete finite-volume
hyperbolic 3-manifold and obtain a contradiction. Usually (but not
always) the contradiction will be that $\Gamma$ is abelian. The
arguments for each group are as follows:

\begin{itemize}
\item
$\langle a,b\ |\ ab^{-1}a^{-2}b^{-1}ab^2\rangle$: by direct calculation,
$[[ab,b],ab]=1$. Therefore $ab$ and $b(ab)^{-1}b^{-1}$ both lie in the
same maximal abelian subgroup of $\Gamma$. By Lemma
\ref{lem:easy_groups} $[ab,b]=1$ and hence $[a,b]=1$.

\item
$\langle a,b\ |\ ab^{-3}ab^2a^2b^2 \rangle$: we get $b^3=(ab^2a)^2$ and
hence $[b^3,ab^2a]=1$. By Lemma \ref{lem:easy_groups},
$[b,ab^2a]=1$. This is equivalent to $[(bab)^2,b^{-1}]=1$; by Lemma
\ref{lem:easy_groups} again we get $[bab,b]=1$ and hence $[a,b]=1$.

\item
$\langle a,b\ |\ ab^2ab^{-1}ab^2a^2b^2ab^{-1} \rangle$: we get
$(ab^2a)^{-2}b=b^{-1}(ab^2a)$. Lemma \ref{lem:group_lem} applies and hence
$[(ab^2a)^{-1},b]=1$. Hence $[ab^2a,b]=1$. Proceed as in the previous
case.

\item
$\langle a,b\ |\ a^3b^2a^3ba^{-2}b \rangle$: we get $a^2=(ba^3b)^2$, so
$[a^2,ba^3b]=1$ and by Lemma \ref{lem:easy_groups} $[a,ba^3b]=1$. This
is equivalent to $[(aba^2)^2,a^{-1}]=1$; by Lemma
\ref{lem:easy_groups} again $[aba^2,a]=1$ and hence $[b,a]=1$.

\item
$\langle a,b,c\ |\ a^2cb^{-1}cb,\ b^2c^2 \rangle$: since $b^2=c^{-2}$,
$[b^2,c^2]=1$ and hence $[b,c]=1$ by Lemma \ref{lem:easy_groups}. The
first relation then simplifies to $a^2c^2=1$, so similarly we get
$[a,c]=1$. By Lemma \ref{lem:easy_groups} $a$ and $b$ must lie in the
unique maximal abelian subgroup containing $c$.

\item
$\langle a,b,c\ |\ bc^{-1}b^{-1}c,\ ac^2ba^{-1}b \rangle$: we have
$[b,c]=1$ immediately, and $a^{-1}ba=(c^2b)^{-1}$, so by Lemma
\ref{lem:easy_groups} $a$ lies in the same abelian subgroup as $b$ and
$c$.

\item
$\langle a,b\ |\ ab^2a^3b^2a^3b^2a^3b^2 \rangle$: we have
$a^2=(a^3b^2)^4$, hence $[a^2,a^3b^2]=1$. By Lemma
\ref{lem:easy_groups} $[a,a^3b^2]=1$ and hence
$[a,b^2]=1$; apply the lemma again to get $[a,b]=1$.

\item
$\langle a,b\ |\ a^2b^2a^2b^2a^3b^2a^2b^2a^2b^2ab^{-1}ab^2 \rangle$:
here we need to alter the group presentation. Let $x=ab^2a$; the presentation
becomes
\[
\langle a,b,x\ |\ x^3ax^3b^{-1},\ ab^2ax^{-1} \rangle .
\]
Then use the first relation to eliminate $b=x^3ax^3$; the presentation
becomes
\[
\langle a,x\ |\ ax^3ax^6ax^3ax^{-1} \rangle .
\]
Now we have $x^6(ax^3a)=(ax^3a)^{-1}x$. By Lemma \ref{lem:group_lem}
$[x^7,ax^3a]=1$; by Lemma \ref{lem:easy_groups} $[x,ax^3a]=1$. This is
equivalent to $[(xax^2)^2,x^{-1}]=1$. By Lemma \ref{lem:easy_groups}
again, $[xax^2,x]=1$ and hence $[a,x]=1$.

\item
$\langle a,b\ |\ a^2ba^2b^{-2}a^2b \rangle$: we get $b^3=(a^2b)^3$, so
$[b^3,a^2b]=1$. Apply Lemma \ref{lem:easy_groups} to get $[b,a^2b]=1$
and hence $[b,a^2]=1$; apply the lemma again to get $[b,a]=1$.

\item
$\langle a,b\ |\ ab^{-3}ab^4a^2b^4 \rangle$: we get $b^3=(ab^4a)^2$, so
$[b^3,ab^4a]=1$ and by Lemma \ref{lem:easy_groups} $[b,ab^4a]=1$. This
is equivalent to $[(bab^3)^2,b^{-1}]=1$; apply the lemma again to get
$[bab^3,b]=1$ and hence $[a,b]=1$.

\item
$\langle a,b,c\ |\ a^2b^{-1}c^{-1}a^{-1}bc,\ b^2c^2 \rangle$: since
$b^2=c^{-2}$, $[b^2,c^2]$=1 and hence $[b,c]=1$ by Lemma
\ref{lem:easy_groups}, i.e. $b$ and $c$ lie in the same maximal
abelian subgroup. Hence either the group has torsion (and is hence not
the fundamental group of a hyperbolic manifold), or else
$b=c^{-1}$. In the latter case, the first relation then implies that
$a$ is trivial, which implies the whole group is abelian.
\end{itemize}

This completes the proof of Proposition \ref{prop:cusped_groups}.\qed

In the closed case, Section 2 produced 69 Dehn fillings for which Snap
failed to find a hyperbolic structure; in addition, there were 3 cases
in which Snap claimed to find a ``nongeometric'' structure on a
non-hyperbolic manifold, as previously discussed.  For these manifolds
we wish to perform a similar analysis to the one just completed for
cusped manifolds.  As in the cusped case, a program written with the
SnapPea kernel was used to re-triangulate these manifolds and attempt
to find manifolds with identical triangulations. This reduces the list
of manifolds from 72 to 24 cases. Fundamental groups were then
calculated for these 24 remaining manifolds. In 15 cases, the
fundamental groups are either cyclic or have a presentation with two
generators and at least one relation of the form $a^nb^m=1$ or
$a^nb^m=b^{-m}a^k$ for some integers $k$, $m$, and $n$. For these
groups, Lemmas \ref{lem:easy_groups} and \ref{lem:group_lem} apply
immediately. The remaining 9 groups require further analysis:

\begin{proposition}
None of the following 9 groups is the fundamental group of a
closed hyperbolic 3-manifold:
\begin{itemize}
\item $\langle a,b\ |\ ababab^{-1}a^2b^{-1},\ 
ab^{-1}ab^2a^{-1}b^2ab^{-1}ab^2 \rangle$
\item $\langle a,b\ |\ ab^3aba^{-2}b,\ 
aba^{-1}bab^{-1}a^{-1}b^{-1} \rangle$
\item $\langle a,b\ |\ ab^{-1}a^{-1}b^2a^{-1}b^{-1}ab,\ 
aba^{-1}b^{-1}a^5b^{-1}a^{-1}b \rangle$
\item $\langle a,b\ |\ a^2b^{-1}ab^{-1}a^2b^2,\
a^2ba^2ba^2b^{-1}a^{-1}b^{-1} \rangle$
\item $\langle a,b\ |\ a^2b^2a^{-1}b^{-1}a^{-1}b^2,\ 
ab^2a^{-1}ba^{-1}b^2ab^{-2}\rangle$
\item $\langle a,b\ |\ abababab^{-1}a^{-1}ba^{-1}b^{-1},\ 
ab^2ab^{-1}a^2b^{-1}\rangle$
\item $\langle a,b\ |\ aba^2bab^{-1}a^2b^{-1},\ 
abab^{-1}a^{-1}ba^{-1}b^{-1} \rangle$
\item $\langle a,b\ |\ ab^{-1}a^{-1}bab^{-1}aba^{-1}b^{-1}ab^2,\ 
aba^{-1}bab^{-1}a^{-1}ba^5ba^{-1}b^{-1} \rangle$
\item $\langle a,b\ |\ a^2ba^{-1}ba^2b^{-1}a^2ba^{-1}b,\ 
ababa^{-1}b^2a^{-1}b \rangle$
\end{itemize}
\label{prop:closed_groups}
\end{proposition}

\noindent\emph{Proof:} As in Proposition \ref{prop:cusped_groups} we
assume that the group is the fundamental group $\Gamma$ of a closed
hyperbolic 3-manifold as use Lemmas \ref{lem:easy_groups} and
\ref{lem:group_lem} to prove the group must be abelian, a
contradiction. The arguments for each group are as follows:

\begin{itemize}
\item
$\langle a,b\ |\ ababab^{-1}a^2b^{-1},\ 
ab^{-1}ab^2a^{-1}b^2ab^{-1}ab^2 \rangle$: from the first relation,
replace $ab^{-1}ab$ with $a^{-1}ba^{-1}b^{-1}a^{-1}$ in the second
relation to get:
\[
\langle a,b\ |\ ababab^{-1}a^2b^{-1},\ 
a^{-1}ba^{-1}b^{-1}a^{-1}ba^{-1}b^2ab^{-1}ab^2 \rangle.
\]
The second relation is now of the form $x^{-2}b^2xb^2$ where
$x=ab^{-1}ab$. Lemma \ref{lem:group_lem} implies that $[x^3,b^2]=1$,
and Lemma \ref{lem:easy_groups} implies that $[x,b]=1$,
i.e. $[ab^{-1}ab,b]=1$, which is equivalent to
$[ab^{-1}ab^{-1},b]=1$. By Lemma \ref{lem:easy_groups} again, we have
$[ab^{-1},b]=1$ and hence $[a,b]=1$.

\item
$\langle a,b\ |\ ab^3aba^{-2}b,\ 
aba^{-1}bab^{-1}a^{-1}b^{-1} \rangle$: the second relation is of the
form $[aba^{-1},b]=1$. By Lemma \ref{lem:easy_groups} we have $[a,b]=1$.

\item
$\langle a,b\ |\ ab^{-1}a^{-1}b^2a^{-1}b^{-1}ab,\ 
aba^{-1}b^{-1}a^5b^{-1}a^{-1}b \rangle$: let $x=[a,b^{-1}]$ and
$y=[a,b]$ to get the following presentation:
\[
\langle a,b,x,y\ |\ xa^{-1}yab,\ ya^4x,\ x^{-1}ab^{-1}a^{-1}b,\ 
y^{-1}aba^{-1}b^{-1} \rangle.
\]
Use the first relation to eliminate $b=a^{-1}y^{-1}ax^{-1}$:
\[
\langle a,x,y\ |\ ya^4x,\ x^{-1}axa^{-1}ya^{-1}y^{-1}ax^{-1},\ 
y^{-2}ax^{-1}a^{-1}xa^{-1}ya \rangle.
\]
Use the new first relation to eliminate $y=x^{-1}a^{-4}$:
\[
\langle a,x\ |\ x^{-1}axa^{-1}x^{-1}a^{-1}xax^{-1},\ 
axa^4xax^{-1}a^{-1}xa^{-1}x^{-1} \rangle.
\]
Use the first relation to replace $axax^{-1}a^{-1}x$ with $xax^{-1}$ in
the second relation to get:
\[
\langle a,x\ |\ x^{-1}axa^{-1}x^{-1}a^{-1}xax^{-1},\ 
axa^3xax^{-1}a^{-1}x^{-1} \rangle.
\]
Use the first relation again to replace $axax^{-1}a^{-1}$ with
$xax^{-2}$ in the second relation:
\[
\langle a,x\ |\ x^{-1}axa^{-1}x^{-1}a^{-1}xax^{-1},\ 
axa^2xax^{-3} \rangle.
\]
Finally make a change of generators by replacing $a$ with $z=ax$ to
get:
\[
\langle x,z\ |\ x^{-1}zxz^{-2}xzx^{-2},\ z^2x^{-1}z^2x^{-4} \rangle.
\]
The second relation now implies that $[x^3,z^2]=1$ by Lemma
\ref{lem:group_lem}, and hence $[x,z]=1$ by Lemma
\ref{lem:easy_groups}.

\item
$\langle a,b\ |\ a^2b^{-1}ab^{-1}a^2b^2,\
a^2ba^2ba^2b^{-1}a^{-1}b^{-1} \rangle$: use the second relation to
replace $ba^2b^{-1}$ with $a^{-2}b^{-1}a^{-2}ba$ in the first relation
to get:
\[
\langle a,b\ |\ a^{-2}b^{-1}a^{-2}ba^2b^{-1}a^2b,\ 
a^2ba^2ba^2b^{-1}a^{-1}b^{-1} \rangle.
\]
The first relation is now of the form $[a^{-2},b^{-1}a^{-2}b]=1$,
hence $b$ and $a^2$ lie in the same maximal abelian subgroup by Lemma
\ref{lem:easy_groups}. Hence by the same lemma $a$ and $b$ lie in the
same maximal abelian subgroup.

\item
$\langle a,b\ |\ a^2b^2a^{-1}b^{-1}a^{-1}b^2,\ 
ab^2a^{-1}ba^{-1}b^2ab^{-2}\rangle$: use the second relation to replace
$ab^2a^{-1}$ with $b^2a^{-1}b^{-2}ab^{-1}$ in the first relation to
get:
\[
\langle a,b\ |\ ab^2a^{-1}b^{-2}ab^{-2}a^{-1}b^2,\ 
ab^2a^{-1}ba^{-1}b^2ab^{-2} \rangle.
\]
The first relation is now of the form $[ab^2a^{-1},b^{-2}]=1$, hence
$a$ and $b^2$ lie in the same maximal abelian subgroup by Lemma
\ref{lem:easy_groups}. Hence by the same lemma $a$ and $b$ lie in the
same maximal abelian subgroup.

\item
$\langle a,b\ |\ abababab^{-1}a^{-1}ba^{-1}b^{-1},\ 
ab^2ab^{-1}a^2b^{-1}\rangle$: the second relation is in the form
$b^3x^2bx^2=1$ where $x=ab^{-1}$. By Lemma \ref{lem:group_lem}, this
implies $[b^2,x^2]=1$, and by Lemma \ref{lem:easy_groups} this implies
that $[b,x]=1$, which implies $[b,a]=1$.

\item
$\langle a,b\ |\ aba^2bab^{-1}a^2b^{-1},\ 
abab^{-1}a^{-1}ba^{-1}b^{-1} \rangle$: the second relation is of the
form $[a,b^{-1}ab]=1$, which implies that $[a,b]=1$ by Lemma
\ref{lem:easy_groups}.

\item
$\langle a,b\ |\ ab^{-1}a^{-1}bab^{-1}aba^{-1}b^{-1}ab^2,\ 
aba^{-1}bab^{-1}a^{-1}ba^5ba^{-1}b^{-1} \rangle$: Let
$x=aba^{-1}b^{-1}$ to get the following presentation:
\[
\langle a,b,x\ |\ x^{-2}axab,\ xbx^{-1}ba^4x,\
x^{-1}aba^{-1}b^{-1}\rangle.
\]
Then use the first relation to eliminate $b=a^{-1}x^{-1}a^{-1}x^2$:
\[
\langle a,x\ |\ xa^{-1}x^{-1}a^{-1}xa^{-1}x^{-1}a^{-1}x^2a^4x,\ 
x^{-2}a^{-1}x^2a^{-1}x^{-2}axa \rangle.
\]
From the second relation, replace $xa^{-1}x^{-1}a^{-1}x^2a$ with
$x^{-1}a^{-1}x^2$ in the first relation to get:
\[
\langle a,x\ |\ xa^{-1}x^{-1}a^{-1}x^{-1}a^{-1}x^2a^3x,\ 
x^{-2}a^{-1}x^2a^{-1}x^{-2}axa \rangle.
\]
Use the second relation again to replace $a^{-1}x^{-1}a^{-1}x^2a$ with
$x^{-2}a^{-1}x^2$ in the first relation to get:
\[
\langle a,x\ |\ xa^{-1}x^{-3}a^{-1}x^2a^2x,\ 
x^{-2}a^{-1}x^2a^{-1}x^{-2}axa \rangle.
\]
Now let $y=x^2$ to get the following presentation:
\[
\langle a,x,y\ |\ ya^{-1}y^{-1}x^{-1}a^{-1}ya^2,\ 
y^{-1}a^{-1}ya^{-1}y^{-1}axa,\ y^{-1}x^2 \rangle.
\]
Use the second relation to eliminate $x=a^{-1}yay^{-1}aya^{-1}$ to get:
\[
\langle a,y\ |\ ya^{-1}y^{-1}ay^{-1}a^{-1}ya,\ 
y^{-1}a^{-1}yay^{-1}aya^{-2}yay^{-1}aya^{-1} \rangle.
\]
The first relation is now of the form $[ya^{-1}y^{-1}a,y^{-1}]=1$,
which implies that $[y,a^{-1}ya]=1$, and hence $[y,a]=1$ by Lemma
\ref{lem:easy_groups}.

\item
$\langle a,b\ |\ a^2ba^{-1}ba^2b^{-1}a^2ba^{-1}b,\ 
ababa^{-1}b^2a^{-1}b \rangle$: this group presentation is just the
first group presentation, with the generators $a$ and $b$ exchanged.
\end{itemize}

This completes the proof of the Proposition \ref{prop:closed_groups}.\qed


\end{document}